\documentclass{article}
\usepackage{amsmath}
\usepackage{amsthm}
\usepackage{amsfonts}
\usepackage{latexsym}
\usepackage{latexsym}
\usepackage[all]{xy}

\newcounter{alphthm}
\newtheorem{thm}{Theorem}[section]
\newtheorem{cor}[thm]{Corollary}
\newtheorem{prop}[thm]{Proposition}

\newtheorem{lem}[thm]{Lemma}

\newcommand{\be}{\begin{equation}}
\newcommand{\ee}{\end{equation}}
\newcommand{\ben}{\begin{enumerate}}
\newcommand{\een}{\end{enumerate}}
\newcommand{\pa}{{\partial}}

\newcommand{\pxi}{{\pa \over \pa x^i}}
\def\beq{\begin{equation}}
\def\eeq{\end{equation}}

\title{\Large  Shen's Processes on Finslerian Connections }
\author{A. Tayebi and  B. Najafi}
\date{}

\begin{document}

\maketitle

\begin{abstract}
In this paper, we discuss the invariant properties of curvatures effected by the Matsumoto's $C$ or $L$-process. We find equivalent conditions on curvatures by comparing the difference between the corresponding curvatures of closely related  connections. As an application, Matsumoto's
$L$-process on Randers manifold is studied. Shen connection can not be obtained by using Matsumoto's processes from other well-known connections. This  leads us to two new processes which we call Shen's $C$ and $L$-processes. We study the invariant properties of curvatures under the Shen's processes.
\end{abstract}

\section{Introduction}
After Einstein's formulation of general relativity, Riemannian geometry became fashionable and one of the
connections, namely Levi-Civita connection, came to forefront. This connection is both torsion-free and metric-compatible. On the other hand, Finsler geometry is a natural extension of Riemannian geometry. Likewise, connections in Finsler geometry can be prescribed on the pulled-back bundle $\pi^* TM$. Examples of such were proposed by Synge, Taylor, Berwald, Cartan and Chern  [1-5],\cite{KT}\cite{TAE}. However,  there are four well-known connections in Finsler geometry which may be considered ``natural" in some sense: Berwald, Cartan, Hashiguchi and Chern connections. Incidentally in the generic Finslerian setting, it is impossible to have a connection on $\pi^* TM$ which is both torsion-free and  compatible with the Riemannian metric  induced by Finsler metric.

In \cite{Ma2}, Matsumoto introduced a satisfactory and truly aesthetical axiomatic description of Cartan's connection in the sixties. After the Cartan connection has been constructed, easy processes, baptized by Matsumoto ``$L$-process" and ``$C$-process" yield the Chern, the Hashiguchi and the Berwald connections.

In this paper, we show that the vv-curvature of connections is invariant under the Matsumoto's  $L$-process. Comparing the corresponding curvatures of related connections obtained by this transformation, we get  equivalent conditions on curvatures.

\begin{thm}\label{thmMatL}
Let $(M,F)$ be a Finsler manifold. Suppose that $\nabla$ and $\widetilde{\nabla}$ are two connections on $M$ and $\widetilde{\nabla}$ is obtained from $\nabla$ by Matsumoto's $L$-process. Then we have the following
\begin{description}
    \item[] (1) Their vv-curvatures coincide.
    \item[] (2) If their hh-curvatures  coincide, then $F$ is a generalized Landsberg metric. Moreover,
 if $M$ is compact, then $F$ reduces to a Landsberg metric.
    \item[] (3) If their hh-curvatures  coincide and $F$ is of non-zero scalar flag curvature, then $F$ is a
Randers metric.
    \item[] (4) Their hv-curvatures  coincide if and only if $F$ is a Landsberg metric.
\end{description}
\end{thm}

It is well known that vanishing hv-curvatures of Cartan and Berwald connections characterize  Landsberg metrics and Berwald metrics, respectively. Shen introduces a new  connection in Finsler geometry, which vanishing  hv-curvature of this  connection  characterizes Riemannian metrics \cite{Sh1}. On the other hand, the Chern, Berwald, and Hashiguchi connections are obtained from Cartan connection by Matsumoto's processes, as depicted in following
\[
\begin{array}{ccccc}
  \textrm{Cartan connection} & \overset{\tiny{C-process}}
  {-\!\!\!-\!\!\!-\!\!\!\longrightarrow}
  & \textrm{Chern connection} &
  \overset{}{} & \textrm{} \\
  | &   & | &  &  \\
  {\tiny\tiny \small L-process}&   & {\tiny {L-process}} &  &  \\
  \downarrow &  & \downarrow &  &  \\
  \textrm{Hashiguchi connection}  & \overset{C-process}{-\!\!\!-\!\!\!-\!\!\!\!
  \longrightarrow} &
  \textrm{Berwald connection} &  &  \\
\end{array}
\]
However, Shen connection can not be constructed by Matsumoto's processes from these well-known connections. Therefore, it is natural to find some kinds of processes on one of these connections, say Chern connection, which yield the Shen connection. Here,  we introduce two new processes on connections called Shen's $C$ and $L$-processes. We show that Shen connection is obtained from Chern connection by Shen's  $C$-process.  Studying curvature tensors of two connections related by this process leads us to the following theorem.
\begin{thm}\label{thmShC}
Let $(M,F)$ be a Finsler manifold. Suppose that $\nabla$ and $\widetilde{\nabla}$ are two connections on $M$ and $\widetilde{\nabla}$  is obtained from $\nabla$ by Shen's $C$-process. Then we have the following
\begin{description}
    \item[](1) If their hh-curvature coincide, then $F$ is a  Landsberg metric.
    \item[](2) Their hv-curvature coincide if and only if $F$ is Riemannian.
    \item[](3) Their vv-curvature coincide.
\end{description}
\end{thm}

Throughout this paper, we set the Cartan connection on Finsler manifolds. The $h$- and $v$- covariant derivatives  are denoted by `` $;$ " and ``, " respectively. Further, we suppose that the horizontal distribution of connections are the same as Cartan connection's horizontal distribution.

\section{Preliminaries}\label{sectionP}
Let $M$ be an n-dimensional $ C^\infty$ manifold. Denote by $T_x M $
the tangent space at $x \in M$, and by $TM=\cup _{x \in M} T_x M $
the tangent bundle of $M$.

A Finsler metric on $M$ is a function $ F:TM
\rightarrow [0,\infty)$ which has the following properties: (i) $F$ is $C^\infty$ on $TM_{0}:= TM \setminus \{ 0 \}$; (ii) $F$ is positively 1-homogeneous on the fibers of tangent bundle $TM$,  and (iii) for each $y\in T_xM$, the following quadratic form $g_y$ on
$T_xM$  is positive definite,
\[
g_{y}(u,v):={1 \over 2} \left[  F^2 (y+su+tv)\right]|_{s,t=0}, \ \
u,v\in T_xM.
\]
Let  $x\in M$ and $F_x:=F|_{T_xM}$.  To measure the
non-Euclidean feature of $F_x$, define ${\bf C}_y:T_xM\times T_xM\times
T_xM\rightarrow \mathbb{R}$ by
\[
{\bf C}_{y}(u,v,w):={1 \over 2} \frac{d}{dt}\left[g_{y+tw}(u,v)
\right]|_{t=0}, \ \ u,v,w\in T_xM.
\]
The family $\bf{C}:=\{\bf{C}_y\}_{y\in TM_0}$  is called the Cartan torsion. It is well known that {\bf{C}=0} if and only if $F$ is Riemannian. For  $y \in T_xM_0$, define  mean Cartan torsion ${\bf I}_y$ by ${\bf I}_y(u) := I_i(y)u^i$, where $I_i:=g^{jk}C_{ijk}$ and $u=u^i{{\partial } \over {\partial x^i}}|_x$. By Diecke's Theorem, $F$ is Riemannian if and only if ${\bf I}_y=0$ \cite{De}.

\smallskip

For  $y \in T_xM_0$, define the  Matsumoto torsion ${\bf M}_y:T_xM\otimes T_xM \otimes T_xM \rightarrow \mathbb{R}$ by ${\bf M}_y(u,v,w):=M_{ijk}(y)u^iv^jw^k$ where
\[
M_{ijk}:=C_{ijk} - {1\over n+1} \Big \{ I_i h_{jk} + I_j h_{ik} + I_k h_{ij} \Big \},\label{Matsumoto}
\]
and $ h_{ij} := FF_{y^iy^j} = g_{ij}- {1\over F^2} g_{ip}y^p g_{jq} y^q $. A Finsler metric $F$ is said to be C-reducible if ${\bf M}_y=0$. This quantity is introduced by  Matsumoto \cite{Ma1}. Matsumoto proves that every Randers metric satisfies that ${\bf M}_y=0$. Later on, Matsumoto-H\={o}j\={o} prove that the converse is true too.
\begin{prop}\label{lemMaHo}{\rm (\cite{Ma3}\cite{MH})}
\emph{A Finsler metric $F$ on a manifold of dimension $n\geq 3$ is a Randers metric if and only if\  ${\bf M}_y =0$, $\forall y\in TM_0$.
}\end{prop}

The horizontal covariant derivatives of ${\bf C}$ and $\bf{I}$ along geodesics give rise to  the  Landsberg curvature  ${\bf L}_y:T_xM\times T_xM\times T_xM\rightarrow \mathbb{R}$ and mean Landsberg curvature  ${\bf J}_y:T_xM\rightarrow \mathbb{R}$ defined by
\[
{\bf L}_y(u,v,w):=L_{ijk}(y)u^iv^jw^k,\ \ \textrm{and}\,\,\,\ {\bf J}_y(u): = J_i (y)u^i,
\]
where $L_{ijk}:=C_{ijk|s}y^s$, $J_i:=I_{i|s}y^s$, $u=u^i{{\partial } \over {\partial x^i}}|_x$,  $v=v^i{{\partial }\over {\partial x^i}}|_x$ and $w=w^i{{\partial }\over {\partial x^i}}|_x$. The families ${\bf L}:=\{{\bf L}_y\}_{y\in TM_{0}}$ and ${\bf J}:=\{{\bf J}_y\}_{y\in TM_{0}}$ are called the Landsberg curvature and mean Landsberg curvature. A Finsler metric is called  Landsberg metric and  weakly Landsberg metric if {\bf{L}=0} and ${\bf J}=0$, respectively.

The rate of change of ${\bf L}$ along geodesics is measured  by the  generalized Landsberg curvature  $\bar{\bf L}_y:T_xM\times T_xM\times T_xM\rightarrow \mathbb{R}$  which is defined by $\bar{\bf L}_y(u,v,w):=\bar{L}_{ijk}(y)u^iv^jw^k$, where $\bar{L}_{ijk}:=L_{ijk|s}y^s$.

The geodesics of Finsler metric $F$ are characterized by the following system of second order ordinary differential equations in local coordinates $ \ddot c^i+2G^i(\dot c)=0,$ where
$G^i(x, y):= {1\over 4} g^{il}(x, y) \{ [F^2]_{x^ky^l}y^k-[F]^2_{x^l}\}$. These local functions $G^i$ define a global vector field on $TM_{0}$ as follows
\[
{\bf G} = y^i\pxi - 2 G^i(x, y) {\pa \over \pa y^i}.
\]
For $y \in T_xM_0$, define ${\bf B}_y:T_xM\otimes T_xM \otimes T_xM\rightarrow T_xM$  by
\[
{\bf B}_y(u, v, w):=B^i_{\ jkl}(y)u^jv^kw^l{{\partial } \over {\partial
x^i}}|_x,
\]
where $ B^i_{\ jkl}(y):={{\partial^3 G^i} \over {\partial y^j \partial y^k \partial y^l}}(y)$,  $u=u^i{{\partial } \over {\partial x^i}}|_x$, $v=v^i{{\partial }\over {\partial x^i}}|_x$ and $w=w^i{{\partial } \over {\partial x^i}}|_x$. $\bf B$ is called the Berwald curvature. A Finsler metric is called a Berwald metric   if $\textbf{B}=0$ \cite{Sh2}. It is well known that every Berwald metric is a Landsberg metric.

\bigskip

The notion of Riemann curvature is extended to Finsler metrics. For  $y\in T_xM_{0}$,
the Riemann curvature ${\bf R}_y: T_xM\to T_xM$ is defined by ${\bf R}_y (u) = R^i_k (y) u^k \; \pxi$ where
\[
 R^i_k (y):= 2 {\pa G^i\over \pa x^k}
- {\pa^2 G^i\over \pa x^j \pa y^k} y^j + 2 G^j {\pa^2 G^i\over \pa y^j \pa y^k} - {\pa G^i\over \pa y^j} {\pa G^j \over \pa y^k}.\label{Rik}
\]

Take an arbitrary plane $P\subset T_xM$ (flag) and a non-zero vector $y\in P$ (flag pole), the  flag curvature $K(P, y)$ is defined by
\[
K(P, y) := { g_y({\bf R}_y (v), v) \over g_y(y,y) g_y(v,v)-g_y(v,y)g_y(v, y) }.
\]
We say that  a Finsler metric $F$ is  of scalar flag curvature if for any $y\in T_xM$,
the flag curvature $K=K(x,y)$ is a scalar function on $TM_0$. If  $K$ is constant, then $F$ is said to be of  constant flag  curvature.

\bigskip

Let us consider the pull-back tangent bundle $\pi^*TM$ over $TM_0$ defined by $ \pi^*TM=\left\{(u, v)\in TM_0 \times TM_0 | \pi(u)\\=\pi(v)\right\} $. Take a local coordinate system $(x^i)$ in $M$, the local natural frame $\{{{\partial} \over {\partial x^i}}\}$  of  $T_xM$ determines a local natural frame $\partial_i|_v$ for
$\pi^*_vTM$ the fibers of  $\pi^*TM$, where  ${\partial _i |_v=(v,{{\partial} \over {\partial x^i}}| _x )}$,
and $v=y^i{{\partial}\over {\partial x^i}}|_x\in TM_0$. The fiber $\pi^*_vTM$ is isomorphic to
 $T_{\pi(v)}M$ where $\pi(v)=x$. There is a canonical section $\ell$  of $\pi^*TM$ defined by $\ell_v=(v,v)/F(v)$.

Let $TTM$ be the tangent bundle of $TM$ and  $\rho$ the canonical linear mapping $\rho:TTM_0\rightarrow \pi^*TM$ defined by $\rho(\hat{X})=(z,\pi_{*}(\hat{X}))$ where $\hat{X}\in T_zTM_0$ and $z\in TM_0$. The bundle map $\rho$ satisfies $ \rho ({\partial \over {\partial x^i}})=\partial_i$ and $ \rho({\partial \over {\partial y^i}})=0$.  Let $ V_zTM$ be the set of vertical vectors at $z$, that is, the set of vectors tangent to the fiber through $z$, or
equivalently $V_zTM=ker \rho$,\  called  the vertical space.

Let $\nabla$ be a linear connection on  $\pi^*TM$. Consider the linear mapping $\mu_z:T_zTM_0\rightarrow T_{\pi z}M$ defined  by $\mu_z(\hat{X})=\nabla_{\hat{X}}F\ell$,\ where $\hat{X}\in T_zTM_0$. The connection $\nabla$ is called a Finsler connection if for every $z\in TM_0$, $\mu_z$ defines an isomorphism of $ V_zTM_0$ onto $T_{\pi z}M$. Therefore, the tangent space $TTM_0$ in $z$ is decomposed as $ T_zTM_0=H_zTM\oplus V_zTM$, where $ H_zTM=\ker\mu_z$
is called the  horizontal space defined by $\nabla$.  Indeed, any tangent vector $\hat{X}\in T_zTM_0$ in $z$ decomposes to $ \hat{X}=H\hat{X}+V\hat{X}$ where $ H\hat{X}\in H_zTM$ and $V\hat{X}\in V_zTM $.

The structural equations of the Finsler connection $\nabla$ are
\begin{eqnarray}
\mathcal{T}(\hat{X},\hat{Y})\!\!\!\!&=&\!\!\!\!\ \nabla_{\hat{X}}Y-\nabla_{\hat{Y}}X-\rho[\hat{X},\hat{Y}],\label{torsion}\\
\Omega(\hat{X},\hat{Y})Z\!\!\!\!&=&\!\!\!\!\ \nabla_{\hat{X}}\nabla_{\hat{Y}}Z- \nabla_{\hat{Y}}\nabla_{\hat{X}}Z
-\nabla_{[\hat{X},\hat{Y}]}Z,\label{curvature}
\end{eqnarray}
where $X=\rho(\hat{X})$, $Y=\rho(\hat{Y})$ and $Z=\rho(\hat{Z})$. The tensors $\mathcal{T}$ and $\Omega$ are called
respectively the Torsion and Curvature tensors of $\nabla$. Three curvature tensors are defined by
 $R(X,Y):=\Omega(H\hat{X},H\hat{Y})$, $P(X,\dot{Y}):=\Omega(H\hat{X},V\hat{Y})$ and  $Q(\dot{X},\dot{Y}):=\Omega(V\hat{X},V\hat{Y})$,\  where $\dot{X}=\mu(\hat{X})$ and $\dot{Y}=\mu(\hat{X})$.

Let $\{e_i\}^n _{i=1}$ be a  local orthonormal (with respect to $g$) frame field for the pulled-back bundle $\pi ^* TM$ such   that $e_n=\ell$. Let $\{\omega^i\}^n_{i=1}$ be its dual co-frame field. One readily finds that $ \omega^n:= {\partial{F} \over {\partial   {y^i}}}dx^i=\omega$, which is called Hilbert form, and    $\omega(\ell)=1$. Put   $\nabla e_i = \omega ^{\ j} _i \otimes e_j$ and $ \Omega e_i=2\Omega ^{\ j} _i \otimes e_j$,   where $\{\Omega ^{\ j} _{i}\}$ and $\{\omega ^{\ j} _{i}\}$ are called respectively,  the  curvature forms and  connection forms of $\nabla$ with respect to   $\{e_{i}\}$. By definition   $\rho =\omega^i \otimes e_i$ and $\mu:=\nabla F\ell=F\omega^{n+i}\otimes e_i$, where $\omega^{n+i}:=\omega^{\ i}_n +d(log F)\delta^i _n$.  It is easy to show that $\{\omega ^i, \omega^{n+i} \}^n_{i=1}$ is a  local basis for $T^*( TM_0).$  In a natural coordinate, we can expand connection forms $\omega^{\ j}_ i$ as follows
\[
\omega^{\ j}_ i:=\Gamma^j_{\ ik}dx^k+F^j_{\ ik}dy^k,
\] where
$\nabla_{\frac{\pa}{\pa x^i}}{}^{\pa_j}=\Gamma^k_{\ ij}\pa_k$ and $\nabla_{\frac{\pa}{\pa y^i}}{}^{\pa_j}=F^k_{\ ij}\pa_k$. In the rest of paper, we suppose that all connections satisfy $F^k_{\ ij}y^i=F^k_{\ ij}y^j=0$.

\bigskip

Let $\{\bar e_i, \dot e_i\}^n _{i=1}$ be the local basis for $T(TM_0)$, which is dual to $\{\omega ^i, \omega^{n+i} \}^n  _{i=1}$, i.e., $\bar e_i \in HTM, \dot e_i \in VTM$ such that  $\rho(\bar e_i)=e_i, \mu(\dot e_i)=F e_i$. Then equations (\ref{torsion}) and (\ref{curvature}) are equivalent to
\begin{eqnarray}
&& d\omega^i-\omega^j \wedge\omega^i_j=\frac{1}{2}S^i_{\ kl}\omega^k \wedge \omega^l+T^i_{\ kl}\omega^k \wedge \omega^{n+l},\label{Torsion}\\
&& d\omega ^{\ j} _i -\omega ^{\ k}_i \wedge \omega^{\ j}_ k = \Omega^{\ j}_ i,\label{Curvature}
\end{eqnarray}
where $T^i_{\ kl}:=\omega^i(\mathcal{T}(\dot{e}_l,\bar {e}_k))$  and $S^i_{\ kl}:=\omega^i(\mathcal{T}(\bar {e}_k,\bar{e}_l))$. Since the $\Omega ^{\ i} _j$ are 2-forms on $TM_0$,  they can be expanded as
\begin{equation}\label{RPQ}
\Omega^{\ j}_ i={1 \over 2}R^{\ j} _{i \ kl} \omega ^k \wedge
\omega^l +P^{\ j} _{i \ kl} \omega ^k \wedge \omega^{n+l}+{1 \over
2}Q^{\ j} _{i \ kl} \omega ^{n+k} \wedge \omega^{n+l}.
  \end{equation}
The objects $R$, $P$ and  $Q$ are called, respectively,  the hh-, hv- and vv-curvature  tensors of $\nabla$ with the
components $R(\bar e_k,\bar e_l)e_i =R^{\ j}_{i \ kl}e_j$, $ P(\bar e_k,\dot e_l)e_i=P^{\ j}_{i \ kl} e_j
$ and $Q(\dot e_k,\dot e_l)e_i=Q^{\ j} _{i \ kl} e_j$. By (\ref{RPQ}), we have $R^{\ j} _{i \ kl}=-R^{\ j} _{i \ lk}$ and   $ Q^{\ j} _{i \  lk}=-Q^{\ j} _{i \ kl}$.

{\section{Matsumoto's $C$ and $L$-processes}
Matsumoto introduces two processes in connection theory that by them, one can construct the Berwald, Hashiguchi and Chern connections from Cartan connection \cite{Ma2}. The space of all connections makes an affine space modeled on the space of $(1,2)$-tensors over pulled-back bundle $\pi^*TM$. It means that adding a $(1,2)$-tensor to a connection makes a new connection. A Finsler metric $F$ gives us two natural $(1,2)$-tensors with components $C^i_{\ jk}$ (=$g^{il}C_{ljk}$) and $L^i_{\ jk}$ (=$g^{il}L_{ljk}$). These two $(1,2)$-tensors play key role in Matsumoto's processes, and in what we call Shen's processes, here. The $C$-processes use Cartan tensor, and the $L$-processes use Landsberg tensor.

\bigskip

Let $(M, F)$ be a Finsler manifold. Suppose that $\nabla$ is a connection with connection forms $\omega^i_j$. We define
\begin{equation}\label{6}
\tilde{\omega}^i_j:=\omega^i_j-C^i_{\ jk} \omega^{n+k}.
\end{equation}
Then $\tilde{\omega}^i_j$ are connection forms of a connection $\widetilde{\nabla}$, that is called the connection obtained from $\nabla$ by Matsumoto's $C$-process. Similarly, we define
\begin{equation}\label{7}
\tilde{\omega}^i_j:=\omega^i_j+L^i_{\ jk} \omega^k.
\end{equation}
Then $\tilde{\omega}^i_j$ are connection forms of a connection $\widetilde{\nabla}$, that is called the connection obtained from $\nabla$ by Matsumoto's $L$-process. Chern  and Hashiguchi connections are obtained from Cartan connection by  Matsumoto's $C$-process, and Matsumoto's $L$-process, respectively.

\subsection{Proof of Theorem \ref{thmMatL}}

First, we recall the following well-known result from \cite{Ma1}.

\begin{lem}\label{Lem}
Let $(M,F)$ be a Finsler manifold and the Cartan tensor satisfies  $C_{ijk}=B_ih_{jk}+B_jh_{ik}+B_kh_{ij}$ such
that $y^iB_i=0$. Then $F$ is a Randers metric.
\end{lem}

To prove the Theorem \ref{thmMatL}, we need the following.
\begin{prop}\label{prop}
Let $(M,F)$ be a generalized Landsberg space. Suppose $c(t)$ is a geodesic. Put ${\bf C}(t)={\bf C}_{\dot c}(U(t),V(t), W(t))$ where $U(t),V(t)$ and $W(t)$ are the parallel vector fields along $c$ . Then following equation holds.
\begin{equation}\label{12}
{\bf C}(t)={\bf L}(0)t+{\bf C}(0).
\end{equation}
\end{prop}
\begin{proof}
Let $p$ be an arbitrary point of $M$, $y,u,v,w\in T_pM$ and $c:(-\infty,\infty)\rightarrow M$ is the unit
speed geodesic passing from $p$ and $\frac{dc}{dt}(0)=y$. If $U(t),V(t)$ and $W(t)$ are the parallel vector fields along $c$ with $U(0)=u, V(0)=v$ and $W(0)=w$, we put
\[
{\bf L}(t)={\bf L}_{\dot c}(U(t),V(t), W(t)).
\]
By definition of Landsberg curvature, we have
\begin{equation}\label{13}
{\bf L}(t)={\bf C}^{'}(t).
\end{equation}
Let
\begin{equation}\label{14}
{\bar{\bf L}(t)}=\bar{\bf L}_{\dot c}(U(t),V(t), W(t)).
\end{equation}
From the definition of ${\bar{\bf L}}_y$, we have
\begin{equation}\label{15}
{\bar{\bf L}}(t)={\bf L}^{'}(t).
\end{equation}
Since $F$ is generalized Landsberg metric, then we have
\begin{equation}
{\bf L}^{'}(t)=0,
\end{equation}
which implies that ${\bf L}(t)={\bf L}(0)$. By (\ref{13}), we get the proof.
\end{proof}

\bigskip

\noindent
\textbf{Proof of Theorem \ref{thmMatL}:}  Let $\widetilde{\nabla}$ be obtained from $\nabla$  by Matsumoto's $L$-process
\[
\widetilde{\omega}^i_j=\omega^i_j+L^i_{\ jk} \omega^k.
\]
Taking an exterior differential of the above relation, yields
\begin{equation}
d\widetilde{\omega}^i_j=d\omega^i_j+dL^i_{\ jk} \wedge \omega^k+L^i_{\ jk}d\omega^k.\label{18}
\end{equation}
On the other hand, we know that
 \begin{equation}
    dL^i_{\ jk}+L^s_{\ jk} \omega ^{\ i} _s - L^i_{\ sk}\omega ^{\ s} _j
-L^i_{\ js}\omega ^{\  s} _k =L^i_{\ jk | s}\omega ^s +L^i_{\ jk.s}\omega ^{n+s},\label{Landsberg}
    \end{equation}
where $``|"$ and $``. "$ denote the horizontal and vertical derivative with respect to $\nabla$. Using (\ref{Torsion}), (\ref{Curvature}), (\ref{18}) and (\ref{Landsberg}), we have
\begin{eqnarray}
\nonumber \widetilde{\Omega}^i_j=\Omega^i_j\!\!\!\!&+&\!\!\!\!\ (L^i_{\ jk|s}\omega^s+L^i_{\ jk.s}\omega^{n+s}-L^s_{\ jk} \omega ^{\ i} _s + L^i_{\ sk}\omega ^{\ s} _j+L^i_{\ js}\omega ^{\ s} _k )\wedge \omega^k\\
\nonumber \!\!\!\!&-&\!\!\!\!\ L^i_{\ ju}(\frac{1}{2}S^u_{\ kl}\omega^l +T^u_{\ kl}\omega^{n+l})\wedge \omega^{k}+ L^i_{\ jk}\omega^s\wedge \omega^k_s\\ \!\!\!\!&-&\!\!\!\!\ (\omega^k_j+L^k_{\ ju}\omega^u)\wedge(\omega^i_k+L^i_{\ km}\omega^m)
+\omega^k_j\wedge \omega ^i_k.\label{19}
\end{eqnarray}
Replacing (\ref{RPQ}) in (\ref{19}) yields
\begin{eqnarray}
\widetilde{R}^i_{j\ kl}\!\!\!\!&=&\!\!\!\!\ R^i_{j\ kl}-(L^i_{\ jk|l}-L^i_{\ jl|k})-(L^m_{\ jk}L^i_{\ ml}-L^m_{\ jl}L^i_{\ mk})+L^i_{\ ju}S^u_{\ kl},\label{20}\\
\widetilde{P}^i_{j\ kl}\!\!\!\!&=&\!\!\!\!\ P^i_{j\ kl}-L^i_{\ jk.l}+L^i_{\ ju}T^u_{\ kl},\label{21} \\
\widetilde{Q}^i_{j\ kl} \!\!\!\!&=&\!\!\!\!\ Q^i_{j\ kl}\label{22}.
\end{eqnarray}
Immediately, we have the proof of part 1. It results that if $\widetilde{\nabla}$ is obtained from $\nabla$ by Matsumoto's $L$-process, then $\nabla$ is torsion-free if and only if $\widetilde{\nabla}$ is torsion-free.\\\\
{\it Proof of part 2.} Let $\widetilde{R}=R$. By (\ref{20}) we have
\begin{equation}\label{23}
L^i_{\ jk|l}=L^i_{\ jl|k}-L^m_{\ jk}L^i_{\ ml}+L^m_{\ jl}L^i_{\ mk}+L^i_{\ ju}S^u_{\ kl}.
\end{equation}
Regularity of $\nabla$ results that $y^l_{\ |k}=0$. Therefore, by contracting with $y^l$, we get $L^i_{\ jk|l}y^l=0$. By our assumption on connections in this paper, we see that $L^i_{\ jk|l}y^l=L^i_{\ jk;l}y^l$. Hence $F$ is a generalized Landsberg metric.

Now, suppose that $M$ is a  compact manifold. By Proposition \ref{prop}, we have the following
\[
{\bf C}(t)={\bf L}(0)t+{\bf C}(0).
\]
Since $M$ is compact then the Cartan tensor is bounded. Using $||C||<\infty$, and letting $t\rightarrow +\infty$ or $t\rightarrow -\infty$, we get $L(0)=L(u,v,w)=0$. It means that $F$ is a Landsberg metric.\\\\
{\it Proof of part 3.} From \cite{Sh2}, for Finsler manifolds of scalar flag curvature we have
\[
L_{ijk|m}y^m=\frac{-F^2}{3}\{K_{.i}h_{jk}+K_{.j}h_{ik}+K_{.k}h_{ij}+3KC_{ijk}\}.
\]
By part 2, $F$ is a generalized Landsberg metric. Then we get
\[
C_{ijk}=\frac{-1}{3K}\{K_{.i}h_{jk}+K_{.j}h_{ik}+K_{.k}h_{ij}\}.
\]
By Lemma  \ref{Lem}, it results that $F$ is a $C$-reducible metric, and by Proposition \ref{lemMaHo},  $F$ is a
Randers metric.\\\\
{\it Proof of part 4.} Suppose that $\widetilde{P}=P$. By (\ref{21}) we have $L^i_{\ jk.l}=L^i_{\ ju}T^u_{\ kl}$. Contracting with $y^k$ yields $L^i_{\ jk}=0$, since Landsberg tensor is positively homogeneous of degree zero and $T^u_{\ kl}y^k=0$. This completes the proof.
\qed

\begin{cor}
Let $\widetilde{\nabla}$  be obtained from $\nabla$ by Matsumoto's $L$-process. Then the hv-curvature of them under $L$-process is invariant if and only if $\nabla$ coincides with $\widetilde{\nabla}$.
\end{cor}

It is obvious that any Landsberg metric is generalized Landsberg metric but the converse is still an open problem.
Following corollary throws a light into this problem.

\begin{cor}\label{cor}
Let $\widetilde{\nabla}$  be obtained from $\nabla$ by Matsumoto's $L$-process. Suppose that their Riemannian curvature coincide. If their hv-curvature are not equal, then $F$ is a generalized Landsberg metric
which is not Landsbergian.
\end{cor}

\bigskip

Now, we consider Matsumoto's $C$-process. By the same argument and technique used in the proof of Theorem \ref{thmMatL}, one can obtain the following theorem.

\begin{thm}\label{thmMatC}
Let $(M,F)$ be a Finsler manifold. Suppose that $\nabla$ and $\widetilde{\nabla}$ are two connections on $M$. Suppose $\widetilde{\nabla}$ is obtained from $\nabla$ by Matsumoto's $C$-process. Then we have the following
\begin{eqnarray}
\widetilde{R}^i_{j\ kl}\!\!\!\!&=&\!\!\!\!\ R^i_{j\ kl}-C^i_{\ ju}R^u_{n\ lk}\label{27},\\
\widetilde{P}^i_{j\ kl}\!\!\!\!&=&\!\!\!\!\ P^i_{j\ kl}-C^i_{\ jl|k}-C^i_{\ ju}P^u_{n\ kl}\label{28}, \\
\widetilde{Q}^i_{j\ kl}\!\!\!\!&=&\!\!\!\!\ Q^i_{j\ kl}+(C^i_{jk.l}-C^i_{jl.k})+(C^u_{jl}C^i_{uk}-C^u_{jk}C^i_{ul})-C^i_{ju}Q^u_{n\ kl}\label{29}.
\end{eqnarray}
\end{thm}

\subsection{Matsumoto's $L$-process on Randers Manifolds}
An $(\alpha, \beta)$-metric is a scalar function on $TM$ defined by
\[
F:=\phi(\frac{\beta}{\alpha})\alpha,\ \ \ s=\beta/\alpha,
\]
where  $\phi=\phi(s)$ is a $C^\infty$ on $(-b_0, b_0)$ with certain regularity, $\alpha=\sqrt{a_{ij}(x)y^iy^j}$ is a Riemannian metric and $\beta =b_i(x)y^i$ is a 1-form on a manifold $M$. Randers metrics are special $(\alpha, \beta)$-metrics which are closely related to Riemannian metrics  defined by $\phi=1+s$, i.e., $F = \alpha +\beta$ and have important applications both in mathematics and physics \cite{Ra}.

In this section, we study the Matsumoto's $L$-process on Randers manifolds equipped with connections whose hh-torsion vanish. We show that the Riemannian curvature of these  connections is invariant under Matsumoto's $L$-process on a Randers manifold $(M, F)$, if and only if $F$ is a Berwald metric. To prove this result, we need the following.

\begin{lem}\label{LemQ}
Let $(M, F)$ be a Finsler manifold and  $\nabla$  be a connection on $M$ satisfying $g_{ij|k}=0$. Suppose that $\widetilde{\nabla}$ is obtained from $\nabla$ by Matsumoto's $L$-process. Then $R=\tilde{R}$ if and only if the following equations hold
\begin{eqnarray}
&&L_{isk}L^s_{\ jl}-L_{isl}L^s_{\ jk}=0,\label{GL4} \\
&&L_{ijl|k}-L_{ijk|l}=0.\label{GL5}
\end{eqnarray}
\end{lem}
\begin{proof} Fix $k$ and $l$ and  put
\[
Q_{ij}:=L_{ijl|k}-L_{ijk|l}+L_{isk}L^s_{\ jl}-L_{isl}L^s_{\ jk}.
\]
One can write
\[
Q_{ij}:=Q^s_{ij}+Q^a_{ij},
\]
where
\[
Q^s_{ij}:=\frac{1}{2}(Q_{ij}+Q_{ji}),\ \ \textrm{and} \ \ Q^a_{ij}:=\frac{1}{2}(Q_{ij}-Q_{ji}).
\]
It is easy to see  that $Q_{ij}=0$ if and only if $Q^s_{ij}=0$ and $Q^a_{ij}=0$. On the other hand, we have
\begin{eqnarray*}
Q_{ji}\!\!\!\!&=&\!\!\!\!\ L_{jil|k}-L_{jik|l}+L_{jsk}L^s_{\ il}-L_{\ jsl}L^s_{\ ik}\\
\!\!\!\!&=&\!\!\!\!\ L_{ijl|k}-L_{ijk|l}+L^s_{\ jk}L_{sil}-L^s_{\ jl}L_{sik}.
\end{eqnarray*}
Hence
\[
Q^s_{ij}=L_{jil|k}-L_{jik|l},
\]
and consequently
\[
Q^a_{ij}=L_{isk}L^s_{\ jl}-L_{isl}L^s_{\ jk}.
\]
This proves the Lemma.
\end{proof}

Now we are ready to prove the mentioned fact.

\begin{thm}\label{MPR}
Let $F=\alpha+\beta$ be a Randers metric on a  manifold $M$ of dimensional $n\geq3$. Suppose that $\nabla$ has vanishing hh-torsion and $g_{ij|k}=0$. Let $\widetilde{\nabla}$ be obtained from $\nabla$ by Matsumoto's $L$-process. Then their hh-curvatures are the same if and only if $F$ is a Berwald metric.
\end{thm}
\begin{proof} Using the assumptions $S^i_{kl}=0$ and $R=\tilde{R}$ in (\ref{20}) imply that
\begin{equation}\label{GL2}
L^i_{\ jl|k}-L^i_{\ jk|l}+L^i_{\ sk}L^s_{\ jl}-L^i_{\ sl}L^s_{\ jk}=0.
\end{equation}
Using the assumption $g_{ij|k}=0$, and lowering indices by $g_{ij}$ imply that (\ref{GL2}) is equivalent to the following
\begin{equation}\label{GL3}
L_{ijl|k}-L_{ijk|l}+L_{isk}L^s_{\ jl}-L_{isl}L^s_{\ jk}=0.
\end{equation}
By Lemma \ref{LemQ}, we have
\begin{eqnarray}
&&L_{isk}L^s_{\ jl}-L_{isl}L^s_{\ jk}=0,\label{GL4}\\
&&L_{ijl|k}-L_{ijk|l}=0.
\end{eqnarray}

A direct computation yields
\begin{eqnarray}
&&h^s_{i}J_s=J_i,\label{GL6}\\
&&h^s_{i}h_{js}=h_{ij},\label{GL6.5}\\
&&g^{ij}h_{ij}=n-1.\label{GL7.5}
\end{eqnarray}
Since $F$ is a Randers metric, then it is C-reducible, i.e.,
\begin{equation}
C_{ijk}={\frac{1}{1+n}}\{h_{ij}I_k+h_{jk}I_i+h_{ki}I_j\},
\end{equation}
Taking a horizontal covariant derivative from above relation, we get
\begin{equation}\label{GL3}
L_{ijk}={\frac{1}{1+n}}\{h_{ij}J_k+h_{jk}J_i+h_{ki}J_j\}.
\end{equation}
Substituting  (\ref{GL3}) into (\ref{GL4}),  one can obtain
\begin{eqnarray}
\{h_{jl}h_{ki}-h_{jk}h_{li}\}J^sJ_s+\{h_{jl}J_k-h_{jk}J_l\}J_i+\{h_{ki}J_l-h_{li}J_k\}J_j=0.\label{GL8}
\end{eqnarray}
Contracting (\ref{GL8}) with $g^{il}g^{jk}$ and using the relations (\ref{GL6}), (\ref{GL6.5}) and (\ref{GL7.5}), we conclude that
\begin{equation}
(n+1)(n-2)J^sJ_s=0.\label{GL10}
\end{equation}
Since $F$ is positive definite and $n>2$, then we have
\begin{equation}
J_s=0.\label{GL11}
\end{equation}
By (\ref{GL3}) and (\ref{GL11}) we conclude that $F$ is a Landsberg metric. It is proved that $F=\alpha+\beta$ is a Landsberg metric if and only if  $F$ is a Berwald metric \cite{Ma1}. This completes  the proof.
\end{proof}

\section{Shen's $C$ and $L$-processes}
Recently, Shen introduced a new torsion-free and almost
metric-compatible connection and proved that hv-curvature of his connection vanishes if and only if the Finsler structure is Riemannian \cite{Sh1}. However, the hv-curvature tensor of the Berwald,  Cartan, Hashiguchi  or the Chern connections does not characterize Riemannian structures.   Shen connection can not be constructed by Matsumoto's processes from Cartan or Chern connection. Therefore it is natural to find  a kind of process on Chern connection which yields the Shen connection. This problem leads us to find two new processes which we call them the Shen's $C$ and $L$-processes.

\bigskip

Let $(M, F)$ be a Finsler manifold. Suppose that $\nabla$ is a connection with connection forms $\omega^i_j$. We define
\begin{equation}\label{31}
\tilde{\omega}^i_j:=\omega^i_j-C^i_{\ jk} \omega^k.
\end{equation}
Then $\tilde{\omega}^i_j$ are connection forms of a connection $\widetilde{\nabla}$, that is called the connection obtained from $\nabla$ by Shen's $C$-process. Similarly, we can define
\begin{equation}\label{32}
\tilde{\omega}^i_j:=\omega^i_j-L^i_{\ jk} \omega^{n+k}.
\end{equation}
Then $\tilde{\omega}^i_j$ are connection forms of a connection $\widetilde{\nabla}$, that is called the connection obtained from $\nabla$ by Shen's $L$-process.

\begin{thm} \label{chernshen}
Shen connection is obtained from the Chern connection by Shen's $C$-process.
\end{thm}

\subsection{Proof of Theorem \ref{thmShC}}

Let $\widetilde{\nabla}$  be obtained from $\nabla$ by Shen's $C$-process. Taking exterior differential from (\ref{31}) yields
\begin{equation}\label{33}
d\widetilde{\omega}^i_j=d\omega^i_j-dC^i_{\ jk} \wedge \omega^k-C^i_{\ jk}d\omega^k.
\end{equation}
On the other hand we have
\begin{equation}\label{Cartan}
dC^i_{\ jk}+C^s_{\ jk} \omega ^{\ i} _s - C^i_{\ sk}\omega ^{\ s} _j-C^i_{\ js}\omega ^{\ s} _k =C^i_{\ jk | s}\omega ^s +C^i_{\ jk.s}\omega ^{n+s},
\end{equation}
where $``|"$ and $``. "$ denote the horizontal and vertical derivative with respect to $\nabla$. Substituting (\ref{Cartan}) into (\ref{33}), and using (\ref{Torsion}) and (\ref{Curvature}) we get
\begin{eqnarray}
\nonumber \widetilde{\Omega}^i_j=\Omega^i_j \!\!\!\!&-&\!\!\!\!\ (C^i_{\ jk|s}\omega^s+C^i_{\ jk.s}\omega^{n+s}-C^s_{\ jk} \omega ^{\ i} _s + C^i_{\ sk}\omega ^{\ s} _j+C^i_{\ js}\omega ^{\  s} _k)\wedge \omega^k\\
\nonumber \!\!\!\!&-&\!\!\!\!\ (\omega^k_j-C^k_{\ jm}\omega^m)\wedge(\omega^i_k-C^i_{\ kl}\omega^l)+\omega^k_j\wedge \omega^i_k-C^i_{\ jk}\omega^s\wedge \omega^k_s\\
\!\!\!\!&+&\!\!\!\!\ C^i_{\ ju}(\frac{1}{2}S^u_{\ kl}\omega^l+T^u_{\ kl}\omega^{n+l})\wedge\omega^k.\label{34}
\end{eqnarray}
Now by decomposing $\widetilde{\Omega}^i_j$ and $\Omega^i_j$ as in (\ref{RPQ}), one can obtain
\begin{eqnarray}
\widetilde{R}^i_{j\ kl}\!\!\!\!&=&\!\!\!\!\ R^i_{j\ kl}+(C^m_{\ jk}C^i_{\ ml}-C^m_{\
jl}C^i_{\ mk})+(C^i_{\ jk|l}-C^i_{\ jl|k})-C^i_{\ ju}S^u_{\ kl}, \label{35}\\
\widetilde{P}^i_{j\ kl}\!\!\!\!&=&\!\!\!\!\ P^i_{j\ kl}-C^i_{\ jk.l}-C^i_{\ jm}T^m_{\ kl},\label{36} \\
\widetilde{Q}^i_{j\ kl}\!\!\!\!&=&\!\!\!\!\ Q^i_{j\ kl}.\label{37}
\end{eqnarray}
{\it{Proof of part 1.}}
Suppose  $R=\widetilde{R}$. Then from (\ref{35}) we have
\begin{equation}\label{38}
C^m_{\ jk}C^i_{\ ml}-C^m_{\ jl}C^i_{\ mk}=C^i_{\ jl|k}-C^i_{\ jk|l}+C^i_{\ ju}S^u_{\ kl}.
\end{equation}
Contracting with $y^l$ yields $L^i_{\ jk}=0$. It means that $F$ is a Landsberg metric.\\\\
{\it{Proof of part 2.}} Suppose that $P=\widetilde{P}$. Then from (\ref{36}) we conclude that
\begin{equation}\label{39}
C^i_{\ jk.l}+C^i_{\ jm}T^m_{\ kl}=0.
\end{equation}
Using the positively homogeneities of Cartan tensor, and contracting (\ref{39}) with $y^l$ yield $C^i_{\ jk}=0$. Therefore, by   Deicke's theorem $F$ is Riemannian. \\

Finally from (\ref{37}), we see that  their vv-curvatures are the same and $\nabla$ is torsion-free  if and only if $\widetilde{\nabla}$ is torsion-free.
\qed

\bigskip

We have some kind of rigidity on Shen's $C$-process.
\begin{cor}
Let $\widetilde{\nabla}$ be obtained from $\nabla$  by Shen's $C$-process. Then the hv-curvature is invariant under Shen's $C$-process if and only if $\nabla=\widetilde{\nabla}$.
\end{cor}

\bigskip

In continue, we study the Shen's  $L$-process and get the following result.
\begin{thm}\label{thmShL}
Let $(M,F)$ be a Finsler manifold. Suppose that $\nabla$ and $\widetilde{\nabla}$
be two connections on $M$ and $\nabla$ is obtained from $\widetilde{\nabla}$ by
Shen's $L$-process. Then we have the following
\begin{eqnarray}
\widetilde{R}^i_{j\ kl}\!\!\!\!&=&\!\!\!\!\ R^i_{j\ kl}-L^i_{\ ju}R^u_{n\ lk}, \label{r4}\\
\widetilde{P}^i_{j\ kl}\!\!\!\!&=&\!\!\!\!\ P^i_{j\ kl}-L^i_{\ jl|k}-L^i_{\ ju}P^u_{n\ kl},\label{p4} \\
\widetilde{Q}^i_{j\ kl}\!\!\!\!&=&\!\!\!\!\ Q^i_{j\ kl}+(L^i_{\ jk.l}-L^i_{\ jl.k})+(L^u_{\ jl}L^i_{\ uk}-L^u_{\ jk}L^i_{\ ul})-L^i_{\ ju}Q^u_{n\ kl}.\label{q4}
\end{eqnarray}
\end{thm}

\smallskip

\begin{cor}
If torsion-free connection $\nabla$ on the Finsler manifold $(M,F)$ remains torsion-free under the Shen's $L$-process, then $F$ is a Landsberg metric. Hence, Shen's $L$-process acts on the set of all torsion-free connections identically.
\end{cor}

\begin{cor}
Let $\nabla$ be obtained from $\widetilde{\nabla}$ by Shen's $L$-process and $\widetilde{\nabla}$ is not torsion-free. If their hv-curvature are equal to zero, then $F$ is a  generalized Landsberg metric such that is not Landsbergian.
\end{cor}

\subsection{Shen's $C$-process on Berwald Connection}
By Theorem \ref{chernshen}, applying Shen's $C$-process on Chern connection gives Shen connection. It is natural to study effect of Shen's $C$-process on the other well-known connections. Here, we study Shen's $C$-process on Berwald connection.

\begin{thm}
Let $(M,F)$ be a Finsler manifold. Suppose that $\nabla$ is the Berwald connection on $M$ and  $\widetilde{ \nabla}$ is obtained from $\nabla$  by Shen's $C$-process. Then  the hv-curvature of $\widetilde{ \nabla}$ vanishes if and only if $F$ is  Riemannian.
\end{thm}
\begin{proof}
The structure  equation of $\widetilde{ \nabla}$ is given by
\begin{eqnarray}
&& d\omega ^i =\omega ^j \wedge \omega^{\ i} _j,\label{nab1}\\
&& dg_{ij}=g_{kj} \omega ^{\ k} _i +g_{ik } \omega ^{\ k} _j+2\{A_{ijk}-L_{ijk}\} \omega ^k +2A_{ijk} \omega ^{n+k},\label{nab2}
\end{eqnarray}
where $A_{ijk}=FC_{ijk}$. Differentiating (\ref{nab2}) and using (\ref{Curvature}), (\ref{nab1}) and (\ref{nab2})  lead to
\begin{eqnarray}
\nonumber g_{kj} \Omega^{\ k}_i+g_{ik}\Omega^{\
k}_j=\!\!\!\!&-&\!\!\!\!\ 2A_{ijk}\Omega^k_n -2A_{ijk | l} \omega ^k \wedge
\omega^l+2A_{ijk.l} \omega ^{n+k} \ \wedge \omega
^{n+l}\\ \nonumber \!\!\!\!&-&\!\!\!\!\ 2\{A_{ijk.l}-A_{ijk | l}\}\omega^k \wedge \omega^{n+l}
\\\!\!\!\!&+&\!\!\!\!\ (L_{ijk | l}\omega^l+L_{ijk . l}\omega ^{n+l})\wedge \omega ^k.
\end{eqnarray}
Using (\ref{RPQ}), yields
\begin{eqnarray}
&& R_{ijkl}+R_{jikl}=-2A_{ijs}R^{\ s}_{n \ kl},\label{nabR}\\
&& P_{ijkl}+P_{jikl}=-2L_{ijk.l}+2\{A_{ijk.l}-A_{ijl| k}\}-2A_{ijs}P^{\ s}_{n \ kl},\label{nabP}\\
&&A_{ijk.l}=A_{ijl.k}.
\end{eqnarray}
Permuting $i,j,k$ in (\ref{nabP}) yields
\begin{eqnarray}
\nonumber P_{ijkl}=-L_{ijk.l}\!\!\!\!&+&\!\!\!\!\ A_{ijk.l}-(A_{ijl | k}+A_{jkl | i}-A_{kil | j})\\\!\!\!\!&+&\!\!\!\!\ A_{kis}P^{\ s} _{n \ jl}-A_{jks} P^{\ s} _{n \ il}-A_{ijs} P^{\ s} _{n\ kl}.\label{nabPP}
  \end{eqnarray}
Multiplying (\ref{nabPP}) with $y^i$ and using $P_{njnl}=0$ yield
\begin{equation}
P_{njkl}=-A_{jkl}.\label{nabPPP}
\end{equation}
By (\ref{nabPPP}) we get the proof.
\end{proof}

\subsection{Shen's $C$-process on Cartan Connection}
Here, we study effect of Shen's $C$-process on the Cartan connection.

\begin{thm}
Let $(M,F)$ be a Finsler manifold. Suppose that $\nabla$ is the Cartan connection on $M$ and  $\widetilde{\nabla}$ is obtained from $\nabla$  by Shen's $C$-process. Then we get
\begin{description}
\item[](1) If hh-curvature of\ \ $\widetilde{\nabla}$ vanishes then $F$ is a Landsberg metric.
\item[](2) The hv-curvature of\ \ $\widetilde{\nabla}$ vanishes if and only if $F$ is   Riemannian.
\end{description}
\end{thm}
\begin{proof}
The structure  equation of $\widetilde{ \nabla}$ is given by
\begin{eqnarray}
&& d\omega ^i =\omega ^j \wedge \omega^{\ i} _j-A^i_{\ kl}\omega^k\wedge \omega^{n+l},\label{nabb1}\\
&& dg_{ij}=g_{kj} \omega ^{\ k} _i +g_{ik } \omega ^{\ k} _j+2A_{ijk}\omega ^k.\label{nabb2}
\end{eqnarray}
 Differentiating (\ref{nabb2}) and using (\ref{Curvature}), (\ref{nabb1}) and (\ref{nabb2})  leads to
 \begin{eqnarray}
\nonumber g_{kj} \Omega^{\ k}_i+g_{ik}\Omega^{\
k}_j=\!\!\!\!&-&\!\!\!\!\  2(A_{ijk |s} \omega ^s+2A_{ijk.s} \omega ^{n+s})\wedge \omega^k\\ \!\!\!\!&-&\!\!\!\!\ 2A_{ijs}A^s_{\ kl}\ \omega^k \wedge \omega^{n+l}.
\end{eqnarray}
Using (\ref{RPQ}), yields
\begin{eqnarray}
&& R_{ijkl}+R_{jikl}=2(A_{ijk|l}-A_{ijl|k}),\label{nabbR}\\
&& P_{ijkl}+P_{jikl}=2(A_{ijk.l}-A_{ijs}A^s_{\ kl}),\label{nabbP}\\
&& Q_{ijkl}+Q_{jikl}=0.
\end{eqnarray}
If the hh-curvature of\ \ $\widetilde{\nabla}$ vanishes, then by (\ref{nabbR}) we have
$A_{ijk|l}=A_{ijl|k}$ which implies that $F$ is a Landsberg metric.

Now let the hv-curvature of\ \ $\widetilde{\nabla}$ vanishes. By (\ref{nabbP}), we get
\begin{eqnarray}
A_{ijk.l}=A_{ijs}A^s_{\ kl}
\end{eqnarray}
Contracting with $y^l$ yields that $F$ is Riemannian.
\end{proof}

\noindent
Akbar Tayebi\\
Faculty  of Science, Department of Mathematics\\
University of Qom\\
Qom. Iran\\
Email:\ akbar.tayebi@gmail.com
\bigskip

\noindent
Behzad Najafi\\
Faculty  of Science, Department of Mathematics\\
Shahed University\\
Tehran. Iran\\
najafi@shahed.ac.ir

\end{document}